\documentclass[reqno,a4paper, 11pt]{amsart}

\usepackage[a4paper=true,pdfpagelabels]{hyperref}
\usepackage{graphicx}

\usepackage[ansinew]{inputenc}
\usepackage{amsfonts,epsfig}
\usepackage{latexsym}
\usepackage{mathabx}
\usepackage{amsmath}
\usepackage{amssymb}

\usepackage{color}   

\newtheorem{theorem}{Theorem}

\newtheorem{corollary}[theorem]{Corollary}

\theoremstyle{definition}

\theoremstyle{remark}

\numberwithin{equation}{section}

\newcommand{\D}{\mathbb{D}}

\newcommand{\SSS}{\mathcal{S}}

\newcommand{\C}{\mathbb{C}}

\newcommand{\e}{\varepsilon}

\renewcommand{\phi}{\varphi}

           \def\e{\varepsilon}
     \def\om{\omega}      
       \def\t{\theta}       
         \def\r{\rho}

\def\U{{\mathcal U}}

\renewcommand{\H}{\mathcal{H}}

\begin{document}

\title[Integral means of derivatives of univalent functions in Hardy spaces]{Integral means of derivatives of univalent functions in Hardy spaces}


\author{Fernando Pérez-González}
\address{Universidad de La Laguna, P.O. Box 456, 38200 La Laguna, Tenerife, Spain}
\email{fernando.perez.gonzalez@ull.es}

\author{Jouni R\"atty\"a}
\address{University of Eastern Finland, P.O. Box 111, 80101 Joensuu, Finland}
\email{jouni.rattya@uef.fi}

\author{Toni Vesikko}
\address{University of Eastern Finland, P.O. Box 111, 80101 Joensuu, Finland}
\email{tonive@uef.fi}

\thanks{The research of the second author was supported in part by Ministerio de
Econom\'{\i}a y Competitividad, Spain, project PGC2018-096166-B-100; La Junta de Andaluc{\'i}a,
projects FQM210  and UMA18-FEDERJA-002.}

\begin{abstract}
We show that the norm in the Hardy space $H^p$ satisfies
	\begin{equation}\label{absteq}
	\|f\|_{H^p}^p\asymp\int_0^1M_q^p(r,f')(1-r)^{p\left(1-\frac1q\right)}\,dr+|f(0)|^p\tag{\dag}
	\end{equation}
for all univalent functions provided that either $q\ge2$ or $\frac{2p}{2+p}<q<2$. This asymptotic was previously known in the cases $0<p\le q<\infty$ and $\frac{p}{1+p}<q<p<2+\frac{2}{157}$ by results due to Pommerenke (1962), Baernstein, Girela and Peláez (2004) and Gonz\'alez and Pel\'aez (2009). It is also shown that \eqref{absteq} is satisfied for all close-to-convex functions if $1\le q<\infty$. A counterpart of \eqref{absteq} in the setting of weighted Bergman spaces is also briefly discussed.
\end{abstract}

\maketitle


\renewcommand{\thefootnote}{}

\section{Introduction and results}

Let $\H(\D)$ denote the space of analytic functions in the unit disc $\D=\{z\in\C:|z|<1\}$. For $0<p\le\infty$, the Hardy space $H^p$ consists of those $f\in\H(\D)$ such that 
	$$
	\|f\|_{H^p}=\sup_{0<r<1}M_p(r,f)<\infty,
	$$
where
	\begin{equation*}
	\begin{split}
	M_p(r,f)=\left(\frac{1}{2\pi}\int_0^{2\pi}|f(re^{i\theta})|^p\,d\theta\right)^\frac1p, \quad 0<r<1,
	\end{split}
	\end{equation*}
is the $L^p$-mean of the restriction of $f$ to the circle of radius $r$, and $M_{\infty}(r,f)=\max_{|z|=r}|f(z)|$ is the maximum modulus function. The monographs \cite{Duren70} and \cite{Garnett} are excellent sources for the theory of the Hardy spaces. 

An injective function in $\H(\D)$ is called a conformal map or univalent, and the class of all such functions is denoted by $\U$. Let $\SSS$ denote the set of $f\in\U$ normalized such that $f(0)=0$ and $f'(0)=1$. If $f\in\U$, then $(f-f(0))/f'(0)$ belongs to $\SSS$. We refer to \cite{Duren83}, \cite{Pommerenke1} and \cite{Pommerenke2} for the theory of univalent functions. 

In 1927 Prawitz~\cite{Prawitz1927} showed that 
		\begin{equation}\label{prawitz}
		M_p^p(r,f)\le p\int_0^rM_{\infty}^p(t,f)\frac{dt}{t},\quad 0<r<1,\quad 0<p<\infty,\quad f\in\SSS.
		\end{equation}
This combined with the Hardy-Littlewood \cite[p.~411]{H-L1932} inequality 
	\begin{equation}\label{eq:hl}
  \int_0^rM_{\infty}^p(t,f)dt\le\pi rM_p^p(r,f),\quad 0<r<1,\quad 0<p<\infty,\quad f\in\mathcal{H}(\D),
  \end{equation}
the proof of which can be found in \cite[Hilfssatz~1]{Po61/62} and \cite[p.~841]{B-G-P2004}, shows that for each $0<p<\infty$ we have the well-known asymptotic equality
	\begin{equation}\label{H^p-1}
	\|f\|_{H^p}^p\asymp\int_0^1M_{\infty}^p(t,f)\,dt,\quad f\in\U.
	\end{equation}
Therefore the containment of $f\in\U$ in the Hardy space $H^p$ is neatly characterized by the behavior of its maximum modulus. Another well-known characterization is given in terms of the arc-length. Namely, the image of the circle of radius $r$ centered at the origin under $f\in\U$ with $f(0)=0$ is a Jordan curve with zero in its inner domain. The length of this image is $2\pi rM_1(r,f')$, and therefore
    \begin{equation}\label{IneqImage}
    M_{\infty}(r,f)\le\pi rM_1(r,f'), \quad 0<r<1,\quad f\in\U,\quad f(0)=0.
    \end{equation}
This inequality is in a sense sharp for conformal maps as its proof shows, and it is actually valid for all analytic functions. Namely, \eqref{eq:hl} applied to $f'$ yields
		\begin{equation*}
    \begin{split}
    |f(re^{i\t})|
    \le\int_0^r|f'(te^{i\t})|\,dt
		\le\pi rM_1\left(r,f'\right),\quad 0<r<1,\quad f\in\H(\D),\quad f(0)=0.
    \end{split}
    \end{equation*}
This combined with the Prawitz' inequality \eqref{prawitz} shows that
	\begin{equation*}
	M_p^p(r,f)\le p\pi^p\int_0^rM_1^p(t,f')t^{1-p}\,dt,\quad 0<r<1,\quad 0<p<\infty,\quad f\in\SSS.
	\end{equation*}
A kind of converse of this inequality is also valid for some $p$. Indeed, in 1962 Pommerenke \cite[Satz~4]{Po61/62} showed that for $0<p<2$ it holds that 
	$$
	r^p\int_0^r M_1^p(t,f')\,dt\lesssim M_p^p(r,f), \quad 0<r<1,\quad f\in\U,\quad f(0)=0.
	$$
This asymptotic inequality, \eqref{H^p-1} and \eqref{IneqImage} show that for each $0<p<2$ we have
	\begin{equation}\label{H^p-2}
	\|f\|_{H^p}^p\asymp\int_0^1M_{1}^p(t,f')\,dt+|f(0)|^p,\quad f\in\U.
	\end{equation}
In 2009 Gonz\'alez and Pel\'aez~\cite[Theorem~1]{GonzalezPelaez2009} generalized this results to the range $0<p<2+\frac{2}{157}$, showed that it fails for $p\ge\frac{100}{17}\approx5.88$, and also observed by using 1967-results due to Thomas~\cite{Thomas1967} that \eqref{H^p-2} is valid for all $0<p<\infty$ if $\U$ is replaced by its proper subclass of all close-to-convex (univalent) functions~\cite[Proposition~1]{GonzalezPelaez2009}. Recall that $f\in\H(\D)$ is close-to-convex if there exists a convex function $g$ such that the real part of the quotient $f'/g'$ is strictly positive on $\D$. The class of close-to-convex functions $f$ normalized such that $f(0)=0$ and $f'(0)=1$ is denoted by $K$ and it was introduced by Kaplan in 1952,~see \cite[Chapter~2]{Duren83} and \cite[Chapter~2]{Pommerenke1} for further information. At this point we only mention that an important subclass of close-to-convex functions is the class of starlike functions. Starlike functions are conformal maps which map $\D$ onto a domain starlike with respect to the origin.

Integral means of derivatives different from $M_1(r,f')$ appearing in \eqref{H^p-2} can also be used to characterize univalent functions in $H^p$. Namely, by combining the 2004-result by Baernstein, Girela and Peláez~\cite[Theorem~1]{B-G-P2004} and \cite[Theorem~2]{GonzalezPelaez2009} due to Gonz\'alez and Pel\'aez we deduce 
	\begin{equation}\label{The-equation-1}
	\|f\|_{H^p}^p\asymp\int_0^1M_q^p(r,f')(1-r)^{p\left(1-\frac1q\right)}\,dr+|f(0)|^p,\quad f\in\U,
	\end{equation}
if either $0<p\le q<\infty$ or $\frac{p}{1+p}<q<p<2+\frac{2}{157}$. The main result of this note shows that these hypotheses can be significantly relaxed in a certain sense.

\begin{theorem}\label{Thm:main1}
Let $0<p,q<\infty$ such that either $\frac{2p}{2+p}<q<2$ or $q\ge2$. Then \eqref{The-equation-1} is valid. Moreover, if $0<p<\infty$ and $1\le q<\infty$, then \eqref{The-equation-1} is valid for all close-to-convex functions $f$.
\end{theorem} 

On one hand, Theorem~\ref{Thm:main1} shows that for $q\ge2$ there is no restriction on $p$. On the other hand, $\frac{2p}{2+p}\in(0,2)$ for all $0<p<\infty$, and hence the range $\frac{2p}{2+p}<q<2$ covers many cases previously excluded by the requirement $p<2+\frac{2}{157}$. However, the hypothesis $\frac{2p}{2+p}<q$ is obviously strictly stronger than $\frac{p}{1+p}<q$ for each $0<p<\infty$. The statement on close-to-convex functions is a generalization of \cite[Proposition~1]{GonzalezPelaez2009} concerning the case $q=1$. 

The proof of Theorem~\ref{Thm:main1} occupies most of the remaining part of the paper, and it is given in Section~\ref{sec2}. At this point we only mention that we offer two proofs concerning the case $q\ge2$, and one of them reveals that for $p\ge q$ we have the asymptotic equality 
	$$
	\|f\|_{H^p}^p\asymp\int_0^1\left(\int_{D(0,r)}\Delta|f'|^q(z)\,dA(z)\right)^\frac{p}{q}(1-r)^p\,dr,\quad f\in\U,
	$$
where, as usual, $\Delta$ stands for the Laplacian. We have not found this asymptotic in the existing literature and believe that it is of interest. 

We next shortly discuss an application of Theorem~\ref{Thm:main1} to Bergman spaces. Let $\om:\D\to[0,\infty)$ such that $\om(z)=\om(|z|)$ for all $z\in\D$, and $\int_\D\om(z)\,dA(z)<\infty$, where $dA(z)$ denotes the element of the Lebesgue area measure on $\D$. For $0<p<\infty$ and such an~$\omega$, the weighted Bergman space $A^p_\om$ consists of $f\in\H(\D)$ such that 
	$$
	\|f\|_{A^p_\om}^p=\int_{\D}|f(z)|^p\om(z)\,dA(z)<\infty.
	$$
By combining Prawitz' result \eqref{prawitz} and the Hardy-Littlewood inequality \eqref{eq:hl}, and then integrating over $[0,1)$ with respect to $\om(r)r\,dr$, we obtain throught Fubini's theorem the chain of inequalities
	\begin{equation*}
	2\int_0^1 M_\infty^p(r,f)\left(\int_r^1\om(t)\,dt\right)\,dr 
	\le\|f\|_{A_{\om}^p}^{p}
	\le2\pi p\int_0^1 M_\infty^p(r,f)\left(\int_r^1\om(t)t\,dt\right)\frac{dr}{r},	
	\end{equation*}
valid for all $f\in\U$ with $f(0)=0$. Standard arguments then show that 
	\begin{equation}\label{Eq:thm-Bergman-2}
	\|f\|_{A_{\om}^p}^{p}
	\asymp\int_0^1 M_\infty^p(r,f)\left(\int_r^1\om(t)t\,dt\right)\,dr,\quad f\in\U.
	\end{equation}
We may also transfer \eqref{The-equation-1} to the setting of the weighted Bergman spaces as the following result shows.

\begin{corollary}\label{cor}
Let $0<p,q<\infty$ and let $\om:\D\to[0,\infty)$ such that $\om(z)=\om(|z|)$ for all $z\in\D$. Further, assume that one of the following conditions is satisfied:
\begin{itemize}
\item[\rm(i)] $0<p\le q<\infty$; 
\item[\rm(ii)] $\frac{p}{1+p}<q<p<2+\frac{2}{157}$; 
\item[\rm(iii)] $q\ge2$; 
\item[\rm(iv)] $\frac{2p}{2+p}<q<2$.
\end{itemize}
Then
	\begin{equation}\label{kjdg}
	\|f\|_{A^p_\om}^p\asymp\int_0^1M_q^p(r,f')(1-r)^{p\left(1-\frac1q\right)}\left(\int_r^1\om(t)t\,dt\right)\,dr+|f(0)|^p\int_0^1\om(r)r\,dr
	,\quad f\in\U.
	\end{equation}
\end{corollary}

The natural approach that we adopt to obtain \eqref{kjdg} consists of first applying \eqref{The-equation-1} to the univalent dilatation $f_r(z)=f(rz)$ appearing in the Bergman space norm of $f$, and then changing the order of radial integrations. The problem then no longer involves the weight $\omega$ and the final step is managed by using the fact $|f'(\rho\xi)|\asymp|f'(r\xi)|$ for all $\xi$ on the boundary of~$\D$ and $0\le r\le\rho<1$ such that $1-r\asymp1-\rho$. Corollary~\ref{cor} is proved in Section~\ref{sec3}.

The case $q=p$ of Corollary~\ref{cor} is of special interest. It states that
	\begin{equation}\label{wuyt}
	\|f\|_{A^p_\om}^p\asymp\int_\D|f'(z)|^p(1-|z|)^{p-1}\left(\int_{|z|}^1\om(r)r\,dr\right)\,dA(z)+|f(0)|^p\int_0^1\om(r)r\,dr,\quad f\in\U.
	\end{equation}
It is well known that this asymptotic equality is valid for all $f\in\H(\D)$ if $\om$ is the standard radial weight $(1-|z|^2)^{\alpha}$ with $-1<\alpha<\infty$. These kind of asymptotic equalities are known as Littlewood-Paley formulas. The rough idea behind these asymptotics  is that $f'$ behaves in a somewhat similar way as $f$ divided by the distance from the boundary. However, it is known that all Bergman spaces do not admit this property. Namely, there exist radial weights $\om$ such that 
	\begin{equation}\label{pipeli}
	\|f\|_{A^p_\om}^p\asymp\int_\D|f'(z)|^p(1-|z|)^{p}W(z)\,dA(z)+|f(0)|^p,\quad f\in\H(\D),
	\end{equation}
fails to be true for each non-negative radial function $W$ on $\D$ unless $p=2$~\cite[Proposition~4.3]{P-R2014}. However, it was recently discovered in \cite[Theorem~5]{PR2021} that \eqref{pipeli} with $W=\om$ is valid if and only if $\om$ satisfies a certain two-sided doubling condition which imposes severe restrictions to the growth, the decay and the oscillation of the weight. Nevertheless, the asymptotic \eqref{wuyt} shows that if we restrict our consideration to univalent functions, then a Littlewood-Paley formula exists for all weighted Bergman spaces induced by radial weights.
	
To this end a couple of words about the notation already used. If there exists a constant $C>0$ such that $A(x)\le C B(x)$ for all $x$ in some set $I$, then we write either $A(x)\lesssim B(x),\,x\in I,$ or $B(x)\gtrsim A(x),\,x\in I$, and the notation $A(x)\asymp B(x)$, $x\in I$, stands for $A(x)\lesssim B(x)\lesssim A(x)$ for all $x\in I$.

\section{Proof of Theorem~\ref{Thm:main1}}\label{sec2}

First observe that \cite[Theorems~6~and~7]{Flett1972} imply 
	\begin{equation}\label{Eq:Flett}
	\int_0^1M_\infty^p(r,f)\,dr\asymp\int_0^1M_\infty^p(r,f')(1-r)^p\,dr+|f(0)|^p,\quad f\in\H(\D).
	\end{equation}
Moreover, a careful inspection of the proof of \cite[Theorem~5.9]{Duren70} shows that for $0<\alpha<\beta\le\infty$ there exists a constant $C=C(\alpha,\beta)>0$ such that
    \begin{equation}\label{Eq:Duren1}
    M_\beta(r,g)\le
    CM_\alpha\left(\frac{1+r}{2},g\right)(1-r)^{\frac1\beta-\frac1\alpha},\quad 0\le
    r<1,\quad g\in\H(\D).
    \end{equation}
By combining \eqref{Eq:Flett} and \eqref{Eq:Duren1}, with $\beta=\infty$ and $\alpha=q$, we deduce
	$$
	\int_0^1M_\infty^p(r,f)\,dr
	\lesssim\int_0^1M_q^p\left(r,f'\right)(1-r)^{p(1-\frac{1}{q})}\,dr+|f(0)|^p
	=I_{p,q}(f)+|f(0)|^p,\quad f\in\H(\D).
	$$
This together with \eqref{H^p-1} yields $\|f\|_{H^p}^p\lesssim I_{p,q}(f)+|f(0)|^p$ for all $f\in\U$. Observe that this part of the proof is valid for all $0<p,q<\infty$.

For the converse implication assume first that $0<q<2$, and write $q=\alpha+\beta$, where $0<\alpha,\beta<q$. By \cite[Proposition~8.1]{Pommerenke2}, for each fixed $0<p<\infty$, we have 
	\begin{equation}\label{Eq:laplacian}
	\int_0^{2\pi}\Delta|f|^p(re^{i\theta})\,d\theta\lesssim\frac{M^p_{\infty}(r,f)}{1-r},\quad \frac{1}{2}\leq r<1,\quad f\in\SSS.
	\end{equation}
H\"older's inequality and \eqref{Eq:laplacian} yield
    \begin{equation*}
    \begin{split}
    2\pi M_q^q(r,f')
    &\le
    \left(\int_0^{2\pi}\left|\frac{f'(re^{i\theta})}{f(re^{i\theta})}\right|^2|f(re^{i\theta})|^{\alpha\frac2{q}}\,d\theta\right)^\frac{q}{2}
    \left(\int_0^{2\pi}|f(re^{i\theta})|^{\beta\frac{2}{2-q}}\,d\theta\right)^{\frac{2-q}{2}}\\
    &=\left(\int_0^{2\pi}|f(re^{i\theta})|^{\alpha\frac{2}{q}-2}|f'(re^{i\theta})|^2\,d\theta\right)^\frac{q}{2}
    \left(\int_0^{2\pi}|f(re^{i\theta})|^{\frac{2\beta}{2-q}}\,d\theta\right)^{\frac{2-q}{2}}\\
    &\lesssim\frac{M^{\alpha}_{\infty}(r,f)}
		{(1-r)^{\frac{q}{2}}}\left(\int_0^{2\pi}|f(re^{i\theta})|^{\frac{2\beta}{2-q}}\,d\theta\right)^{\frac{2-q}{2}},\quad \frac12\le r<1,\quad f\in\SSS.
    \end{split}
    \end{equation*}
Another application of H\"older's inequality gives
		\begin{equation}\label{hgfd}
    \begin{split}
    I_{p,q}(f)
		&\lesssim\int_{\frac12}^1M_q^p\left(r,f'\right)(1-r)^{p(1-\frac{1}{q})}\,dr\\
    &\lesssim\int_{\frac12}^1 M_\infty^{\alpha\frac{p}{q}}(r,f)
		\left(\int_0^{2\pi}|f(re^{i\theta})|^{\frac{2\beta}{2-q}}\,d\theta\right)^{\frac{2-q}{2}\frac{p}{q}}
		(1-r)^{p\frac{q-2}{2q}}\,dr\\
		&\lesssim\left(\int_0^1M_\infty^p(r,f)\,dr\right)^\frac{\alpha}{q}
		\left(\int_0^1\left(\int_0^{2\pi}|f(re^{i\theta})|^{\frac{2\beta}{2-q}}\,d\theta\right)^{p\frac{2-q}{2\beta}}
		(1-r)^{-p\frac{2-q}{2\beta}}\,dr\right)^\frac{\beta}{q}\\
		&\asymp\left(\int_0^1M_\infty^p(r,f)\,dr\right)^\frac{\alpha}{q}
		\left(\int_0^1M_s^p(r,f)(1-r)^{-\frac{p}{s}}\,dr\right)^\frac{\beta}{q},
		\end{split}
    \end{equation}
where $s=\frac{2\beta}{2-q}$. Our hypothesis $q>\frac{2p}{2+p}$ allows us to choose $\frac{(2-q)p}{2}<\beta<q$ which guarantees $s>p$. Now \cite[Theorems~6~and~7]{Flett1972} imply
	\begin{equation}\label{iuyt}
	\int_0^1M_s^p(r,f)(1-r)^{-\frac{p}{s}}\,dr
	\asymp\int_0^1M_s^p(r,f')(1-r)^{p\left(1-\frac{1}{s}\right)}\,dr+|f(0)|^p,\quad f\in\H(\D),
	\end{equation}
where the right-hand side is comparable to $\|f\|_{H^p}^p$ for all $f\in\U$ by \eqref{The-equation-1}. Therefore \eqref{hgfd}, \eqref{iuyt}, \eqref{The-equation-1} and \eqref{H^p-1} give $I_{p,q}(f)\lesssim\|f\|_{H^p}^p$ for all $f\in\SSS$. An application of this to $(f-f(0))/f'(0)$ gives the assertion for $f\in\U$.

Assume now that $2\le q<\infty$. We offer two different proofs of which the first one is valid for all $0<p<\infty$ and the second one only for $q\le p$. The importance of the second proof lies in the fact that it allows us to characterize univalent functions in $H^p$ in terms of the Laplacian of $|f'|^q$. The first proof is pretty straightforward and reads as follows. It is well-known~\cite[Chapter~5]{Duren70} that for each fixed $0<p\le\infty$ we have 
	\begin{equation}\label{Eq:maxmod}
	M_p(r,f')\lesssim\frac{M_p\left(\frac{1+r}{2},f\right)}{1-r},\quad 0<r<1,\quad f\in\H(\D).
	\end{equation}
	By applying \eqref{Eq:maxmod}, with $p=\infty$, and \eqref{Eq:laplacian} we deduce
	\begin{equation*}
	\begin{split}
	I_{p,q}(f)
	&\asymp\int_{\frac12}^1M_q^p\left(r,f'\right)(1-r)^{p(1-\frac{1}{q})}\,dr\\
	&\lesssim\int_{\frac12}^1\left(M_\infty(r,f')\right)^{(q-2)\frac{p}{q}}
	\left(\int_0^{2\pi}|f'(re^{i\theta})|^2\,d\theta\right)^\frac{p}{q}(1-r)^{p(1-\frac1q)}\,dr\\
	&\lesssim\int_0^1\frac{\left(M_\infty\left(\frac{1+r}{2},f\right)\right)^{(q-2)\frac{p}{q}}}{(1-r)^{(q-2)\frac{p}{q}}}
	\left(\frac{M^2_\infty(r,f)}{1-r}\right)^\frac{p}{q}(1-r)^{p(1-\frac1q)}\,dr\\
	&\lesssim\int_0^1M^p_\infty\left(\frac{1+r}{2},f\right)\,dr
	\lesssim\int_0^1M^p_\infty\left(r,f\right)\,dr,\quad f\in\SSS.
	\end{split}
	\end{equation*}
This together with \eqref{H^p-1} yields $I_{p,q}(f)+|f(0)|^p\lesssim\|f\|_{H^p}^p$ for all $f\in\U$.

The first step towards the second proof is to estimate $I_{p,q}(f)$ upwards for all $f\in\H(\D)$, and it is valid on the range $0<\frac{p}{p+1}<q\le p<\infty$, that is, $0<p(1-\frac1q)+1\le p$. An integration by parts and H\"older's inequality show that
	\begin{equation*}
	\begin{split}
	I_{p,q}(f)&\lesssim\int_0^1\frac{\partial}{\partial r}M_q^p(r,f')(1-r)^{p(1-\frac1q)+1}\,dr+|f'(0)|^p\\
	&=\frac{p}{q}\int_0^1M_q^{p-q}(r,f')\left(\frac{\partial}{\partial r}M_q^q(r,f')\right)(1-r)^{p(1-\frac1q)+1}\,dr+|f'(0)|^p\\
	&\lesssim\left(I_{p,q}(f)\right)^\frac{p-q}{p}
	\left(\int_0^1\left(\frac{\partial}{\partial r}M_q^q(r,f')\right)^\frac{p}{q}(1-r)^p\,dr\right)^\frac{q}{p}+|f'(0)|^p,\quad f\in\H(\D),
	\end{split}
	\end{equation*}
and it follows that 
	\begin{equation}\label{ert}
	I_{p,q}(f)\lesssim\int_0^1\left(\frac{\partial}{\partial r}M_q^q(r,f')\right)^\frac{p}{q}(1-r)^p\,dr+|f'(0)|^p,\quad f\in\H(\D),
	\end{equation}
provided $0<p(1-\frac1q)+1<p$. The next step is to estimate $\frac{\partial}{\partial r}M_q^q(r,f')$ for $f\in\U$, and this is done in three separate cases. If $q\ge4$, then \eqref{Eq:maxmod}, with $p=\infty$, yields
	\begin{equation*}
	\begin{split}
	2\pi r\frac{\partial}{\partial r}M_q^q(r,f')
	&=q^2\int_{D(0,r)}|f'(z)|^{q-2}|f''(z)|^2\,dA(z)\\
	&\le q^2M^2_\infty(r,f'')M_\infty^{q-4}(r,f')\int_{D(0,r)}|f'(z)|^2\,dA(z)\\
	&\lesssim\frac{M_\infty^{q-2}\left(\frac{1+r}{2},f'\right)}{(1-r)^{2}}M^2_\infty(f,r)
	\lesssim\frac{M_\infty^{q}\left(\frac{3+r}{4},f\right)}{(1-r)^{q}},\quad 0<r<1,\quad f\in\U.
	\end{split}
	\end{equation*}
If $2<q<4$, then H\"older's inequality and \eqref{Eq:maxmod}, first with $p=\frac{4}{4-q}$ and then for $p=\infty$, yield
	\begin{equation*}
	\begin{split}
	2\pi r\frac{\partial}{\partial r}M_q^q(r,f')
	&\le q^2\left(\int_{D(0,r)}|f'(z)|^2\,dA(z)\right)^\frac{q-2}{2}
	\left(\int_{D(0,r)}|f''(z)|^\frac{4}{4-q}\,dA(z)\right)^\frac{4-q}{2}\\
	&\lesssim M^{q-2}_\infty(r,f)\left(\int_0^r\frac{M_{\frac{4}{4-q}}^{\frac{4}{4-q}}\left(\frac{1+s}{2},f'\right)}{(1-s)^{\frac{4}{4-q}}}\,ds\right)^{\frac{4-q}{2}}\\
	&\lesssim\frac{M^{q-2}_\infty(r,f)}{(1-r)^{2}}
	\left(M_\infty\left(\frac{1+r}{2},f'\right)\right)^{\frac{4-q}{2}\left(\frac4{4-q}-2\right)}M^{4-q}_\infty\left(\frac{1+r}{2},f\right)\\
	&\le\frac{M^{2}_\infty\left(\frac{1+r}{2},f\right)}{(1-r)^{2}}M^{q-2}_\infty\left(\frac{1+r}{2},f'\right)\\
	&\lesssim\frac{M^{q}_\infty\left(\frac{3+r}{4},f\right)}{(1-r)^{q}},\quad 0<r<1,\quad f\in\U.
	\end{split}
	\end{equation*}
In the case $q=2$, \eqref{Eq:maxmod}, with $p=2$, gives 
	\begin{equation*}
	\begin{split}
	2\pi r\frac{\partial}{\partial r}M_2^2(r,f')
	&=4\int_{D(0,r)}|f''(z)|^2\,dA(z)
	\lesssim\int_0^r\frac{M_2^2\left(\frac{1+s}{2},f'\right)}{(1-s)^2}\,ds\\
	&\lesssim\frac{M_\infty^{2}\left(\frac{1+r}{2},f\right)}{(1-r)^{2}},\quad 0<r<1,\quad f\in\U.
	\end{split}
	\end{equation*}
Therefore we have shown that, for each fixed $2\le q<\infty$, we have
	\begin{equation}\label{bgt}
	2\pi r\frac{\partial}{\partial r}M_q^q(r,f')
	=\int_{D(0,r)}\Delta|f'|^q(z)\,dA(z)\lesssim\frac{M_\infty^q\left(\frac{3+r}{4},f\right)}{(1-r)^q},\quad 0<r<1,\quad f\in\U.
	\end{equation}
This estimate together with \eqref{ert} and \eqref{H^p-1} yields
	\begin{equation*}
	\begin{split}
	I_{p,q}(f)
	&\lesssim\int_\frac12^1\left(\frac{\partial}{\partial r}M_q^q(r,f')\right)^\frac{p}{q}(1-r)^p\,dr+|f'(0)|^p
	\lesssim\int_0^1M_\infty^p\left(\frac{3+r}{4},f\right)\,dr+|f'(0)|^p\\
	&\lesssim\int_0^1M_\infty^p\left(r,f\right)\,dr+|f'(0)|^p\lesssim\|f\|_{H^p}^p,\quad f\in\U.
	\end{split}
	\end{equation*}
This finishes the proof of \eqref{The-equation-1} for $q\ge2$. 

It remains to prove the assertion for close-to-convex functions. Since we have already shown that $\|f\|_{H^p}^p\lesssim I_{p,q}(f)+|f(0)|^p$ for all $f\in\U$, and this is valid for all $0<p,q<\infty$, it remains to estimate $I_{p,q}(f)$ upwards to $\|f\|_{H^p}^p$. We claim that for $1\le q<\infty$, $\e>0$ and for all close-to-convex functions $f$ we have
	\begin{equation}\label{Eq:ctc}
	M_q^q(r,f')\lesssim\frac{1}{(1-r)^\e}\int_0^r\frac{M_\infty^q(t,f)}{(1-t)^{q-\e}}\,dt,\quad \frac12\le r<1,
	\end{equation}
the proof of which is postponed for a moment. Let $\frac{p-q}{q}<x<p$, and pick up $\e=\e(p,q)>0$ such that $x<p(1-\frac{\e}{q})$. Then \eqref{Eq:ctc}, H\"older's inequality and Fubini's theorem yield
	\begin{equation*}
	\begin{split}
	I_{p,q}(f)
	&\lesssim\int_\frac12^1\left(\int_0^r\frac{M_\infty^q(t,f)}{(1-t)^{q-\e-\frac{qx}{p}}}(1-t)^{-\frac{qx}{p}}\,dt\right)^{\frac{p}{q}}
	(1-r)^{p(1-\frac{1+\e}q)}\,dr\\
	&\le\int_0^1\frac{M^p_\infty(t,f)}{(1-t)^{p-\frac{\e p}{q}-x}}\left(\int_t^1\left(\int_0^r\frac{ds}{(1-s)^\frac{qx}{p-q}}\right)^{\frac{p-q}{q}}
	(1-r)^{p(1-\frac{1+\e}q)}\,dr\right)dt\\
	&\asymp\int_0^1\frac{M^p_\infty(t,f)}{(1-t)^{p-\frac{\e p}{q}-x}}\left(\int_t^1\frac{dr}{(1-r)^{x+1-p+\frac{p\e}{q}}}\right)dt
	\lesssim\int_0^1M^p_\infty(t,f)\,dt,
	\end{split}
	\end{equation*}
and we are done by \eqref{H^p-1}. 

It remains to prove \eqref{Eq:ctc}. If $2\le q<\infty$, this, with $\e=0$, follows from the estimate \eqref{bgt} by integrating. Namely, \cite[Corollary~1.6]{Pommerenke2} shows that $|f'(\rho\xi)|\asymp|f'(r\xi)|$ for all $\xi$ on the boundary of $\D$ and $0\le r\le\rho<1$ such that $1-r\asymp1-\rho$, provided $f\in\U$, and hence $M_q(r,f')\asymp M_q(\r,f')$. This together with \eqref{bgt} and \eqref{Eq:maxmod} imply
	\begin{equation*}
	\begin{split}
	M_q^q\left(\frac{3+r}{4},f'\right)
	&\asymp M_q^q(r,f')
	\lesssim\int_0^{r}\frac{M_\infty^q\left(\frac{3+t}{4},f\right)}{(1-t)^q}dt+|f'(0)|^q\\
	&\lesssim\int_{\frac34}^{\frac{3+r}{4}}\frac{M_\infty^q\left(t,f\right)}{(1-t)^q}dt+M_q^q\left(\frac12,f\right)
	\lesssim\int_{0}^{\frac{3+r}{4}}\frac{M_\infty^q\left(t,f\right)}{(1-t)^q}dt,\quad 0<r<1.
	\end{split}
	\end{equation*}
Further, by the proofs of \cite[Theorems~2 and 3]{Thomas1967} we have \eqref{Eq:ctc} for $q=1$ with $\e=0$. Since the right-hand side of \eqref{Eq:ctc} is increasing, it remains to consider the case $1<q<2$. To do this, we use ideas from \cite{Thomas1967}. Since $f$ is close-to-convex, Alexander's theorem~\cite[Theorem~2.12]{Duren83} implies that there exists a starlike function $g$ such that $\Re\frac{zf'(z)}{g(z)}>0$ for all $z\in\D$. Write $h(z)=\frac{zf'(z)}{g(z)}$ for all $z\in\D$. Then
	\begin{equation*}
	\begin{split}
	r^q\int_0^{2\pi}|f'(re^{i\theta})|^q\,d\theta
	&=\int_0^{2\pi}|g(re^{i\theta})h(re^{i\theta})|^q\,d\theta\\
	&\lesssim\int_0^{2\pi}\left(\int_0^r|g'(te^{i\theta})h(te^{i\theta})|\,dt\right)^q\,d\theta
	+\int_0^{2\pi}\left(\int_0^r|g(te^{i\theta})h'(te^{i\theta})|\,dt\right)^q\,d\theta\\
	&=I_1(r)+I_2(r).
	\end{split}
	\end{equation*}
Since $g$ is starlike, there exists $\phi$ such that $\Re\phi>0$ and $zg'(z)=g(z)\phi(z)$ for all $z\in\D$. Hence
	\begin{equation*}
	\begin{split}
	I_1(r)&=\int_0^{2\pi}\left(\int_0^r|g'(te^{i\theta})h(te^{i\theta})|\,dt\right)^q\,d\theta
	=\int_0^{2\pi}\left(\int_0^r|f'(te^{i\theta})\phi(te^{i\theta})|\,dt\right)^q\,d\theta.
	\end{split}
	\end{equation*}
Let $x=\frac{\e+q-1}{q}>0$. Observe that $|\widehat{\phi}(n)|\lesssim1$ because $\Re\phi>0$. This together with H\"older's inequality, Fubini's theorem, \eqref{Eq:laplacian}, \eqref{Eq:maxmod} and Parseval's identity yields
	\begin{equation*}
	\begin{split}
	I_1(r)
	&=\int_0^{2\pi}\left(\int_0^r|f'(te^{i\theta})\phi(te^{i\theta})|(1-t)^x(1-t)^{-x}\,dt\right)^q\,d\theta\\
	&\lesssim\frac{1}{(1-r)^\e}\int_0^r\left(\int_0^{2\pi}|f'(te^{i\theta})\phi(te^{i\theta})|^q\,d\theta\right)(1-t)^{q-1+\e}\,dt\\
	&\lesssim\frac{1}{(1-r)^\e}\int_0^rM_\infty^{q-1}(t,f')M_\infty^{q-1}(t,\phi)\\
	&\quad\cdot\left(\int_0^{2\pi}|f'(te^{i\theta})|^2\,d\theta\right)^\frac12
	\left(\int_0^{2\pi}|\phi(te^{i\theta})|^2\,d\theta\right)^\frac12(1-t)^{q-1+\e}\,dt\\
	&\lesssim\frac{1}{(1-r)^\e}\int_0^r\frac{M_\infty^{q-1}\left(\frac{1+t}{2},f\right)}{(1-t)^{q-1}}\left(\sum_n|\widehat{\phi}(n)|t^n\right)^{q-1}\\
	&\quad\cdot\left(\frac{M_\infty^2(t,f)}{1-t}\right)^\frac12
	\left(\sum_n|\widehat{\phi}(n)|^2t^{2n}\right)^\frac12(1-t)^{q-1+\e}\,dt\\
	&\lesssim\frac{1}{(1-r)^\e}\int_0^r\frac{M_\infty^{q}\left(\frac{1+t}{2},f\right)}{(1-t)^{q-\e}}\,dt
	\lesssim\frac{1}{(1-r)^\e}\int_{0}^{\frac{1+r}{2}}\frac{M_\infty^{q}\left(t,f\right)}{(1-t)^{q-\e}}\,dt.
	\end{split}
	\end{equation*}
Since $\Re h>0$, there exists an increasing function $\mu$ such that
	$$
	h(z)=\frac1{2\pi}\int_0^{2\pi}\frac{1+ze^{-i\alpha}}{1-ze^{-i\alpha}}d\mu(\alpha),\quad \frac1{2\pi}\int_0^{2\pi}d\mu(\alpha)=1.
	$$
Hence
	$$
	h'(z)=\frac1{\pi}\int_0^{2\pi}\frac{e^{-i\alpha}}{(1-ze^{-i\alpha})^2}d\mu(\alpha).
	$$
By using \cite[Lemma~2]{Thomas1967} and \eqref{Eq:maxmod} it follows that
	\begin{equation*}
	\begin{split}
	I_2(r)
	&\lesssim\frac{1}{(1-r)^\e}\int_0^r\left(\int_0^{2\pi}|g(te^{i\theta})h'(te^{i\theta})|^q\,d\theta\right)(1-t)^{q-1+\e}\,dt\\
	&\le\frac{1}{(1-r)^\e}\int_0^r\left(\int_0^{2\pi}\left(\int_0^{2\pi}\frac{1-t^2}{|1-te^{i\theta}e^{-i\alpha}|^2}\,d\mu(\alpha)\frac{|g(te^{i\theta})|}{1-t^2}\right)^q\,d\theta\right)(1-t)^{q-1+\e}\,dt\\
	&\le\frac{1}{(1-r)^\e}\int_0^r\left(\int_0^{2\pi}\left(\Re h(te^{i\theta})|g(te^{i\theta})|\right)^q\,d\theta\right)(1-t)^{\e-1}\,dt\\
	&\lesssim\frac{1}{(1-r)^\e}\int_0^rM^{q-1}_\infty(t,f')\left(\int_0^{2\pi}\Re h(te^{i\theta})|g(te^{i\theta})|\,d\theta\right)(1-t)^{\e-1}\,dt\\
	&\lesssim\frac{1}{(1-r)^\e}\int_0^r\frac{M^{q-1}_\infty\left(\frac{1+t}{2},f\right)}{(1-t)^{q-1}}M_\infty(r,f)(1-t)^{\e-1}\,dt
	\lesssim\frac{1}{(1-r)^\e}\int_{0}^{\frac{1+r}{2}}\frac{M_\infty^{q}\left(t,f\right)}{(1-t)^{q-\e}}\,dt.
	\end{split}
	\end{equation*}
Therefore we have shown that 
	$$
	r^qM_q^q(r,f')\lesssim I_1(r)+I_2(r)\lesssim\frac{1}{(1-r)^\e}\int_{0}^{\frac{1+r}{2}}\frac{M_\infty^{q}\left(t,f\right)}{(1-t)^{q-\e}}\,dt,\quad 0<r<1,
	$$
provided $1<q<2$. This together with \cite[Corollary~1.6]{Pommerenke2} yields	\eqref{Eq:ctc}.

\section{Proof of Corollary~\ref{cor}}\label{sec3}

If the integral $\int_0^1\om(r)r\,dr$ vanishes or diverges then there is nothing to prove, so assume $\int_0^1\om(r)r\,dr\in(0,\infty)$. Then there exists $R=R(\om)\in(0,1)$ such that $\int_{R}^1\om(r)r\,dr>0$. By \eqref{The-equation-1} and Theorem~\ref{Thm:main1} we have
	\begin{equation*}
	\begin{split}
	\|f\|_{A^p_\om}^p
	&\asymp\int_{R}^1\|f_r\|_{H^p}^p\om(r)r\,dr
	\asymp\int_R^1\left(\int_0^1M_q^p(s,f_r')(1-s)^{p\left(1-\frac1q\right)}r^p\,ds+|f_r(0)|^p\right)\om(r)r\,dr\\
	&=\int_R^1\left(\int_0^rM_q^p(t,f')\left(1-\frac{t}{r}\right)^{p\left(1-\frac1q\right)}r^{p-1}\,dt+|f(0)|^p\right)\om(r)r\,dr\\
	&\asymp\int_R^1\left(\int_0^rM_q^p(t,f')\left(r-t\right)^{p\left(1-\frac1q\right)}\,dt\right)\om(r)r\,dr+|f(0)|^p\int_0^1\om(r)r\,dr,\quad f\in\U.
	\end{split}
	\end{equation*}
Moreover, Fubini's theorem yields
	\begin{equation*}
	\begin{split}
	&\int_0^1M_q^p(t,f')\left(1-t\right)^{p\left(1-\frac1q\right)}\left(\int_t^1\om(r)r\,dr\right)dt\\
	&=\int_0^1\left(\int_0^rM_q^p(t,f')\left(1-t\right)^{p\left(1-\frac1q\right)}\,dt\right)\om(r)r\,dr\\
	&\asymp\int_R^1\left(\int_0^rM_q^p(t,f')\left(1-t\right)^{p\left(1-\frac1q\right)}\,dt\right)\om(r)r\,dr.
	\end{split}
	\end{equation*}
Therefore it suffices to show that
	\begin{equation}\label{alpha}
	\int_0^rM_q^p(t,f')\left(1-t\right)^{\alpha}\,dt
	\asymp\int_0^rM_q^p(t,f')\left(r-t\right)^{\alpha}\,dt,\quad R\le r<1,
	\end{equation}
where $\alpha=\alpha(p,q)=p\left(1-\frac1q\right)$. To see this, fix $M=M(\om)>0$ such that $M>1/R$. If $t\le\frac{Mr-1}{M-1}\in(0,r)$, then $1-t\le M(r-t)\le M(1-t)$. Moreover, \cite[Corollary~1.6]{Pommerenke2} yields
	\begin{equation*}
	\int_{\frac{Mr-1}{M-1}}^rM_q^p(t,f')\left(1-t\right)^{\alpha}\,dt
	\asymp M_q^p(r,f')(1-r)^{\alpha+1}
	\asymp \int_{\frac{Mr-1}{M-1}}^rM_q^p(t,f')\left(r-t\right)^{\alpha}\,dt,
	\end{equation*}
from which \eqref{alpha} follows.

\end{document}